\documentclass[12pt]{article}
\usepackage{amsfonts}
\usepackage{amsbsy}
\usepackage{latexsym}
\usepackage{graphicx}
\newcommand{\ov}{\overline}

\newcommand{\al}{\alpha}

\newcommand{\de}{\delta}

\newcommand{\pr}{\prime}

\newcommand{\be}{\begin{equation}}
\newcommand{\ee}{\end{equation}}
\newcommand{\bq}{\begin{eqnarray}}
\newcommand{\eq}{\end{eqnarray}}

\newcommand{\ba}{\begin{array}}
\newcommand{\ea}{\end{array}}

\newcommand{\cpo}{{\mathbb {CP}}^1}

\newcommand{\la}{\lambda}
\newcommand{\bt}{\beta}

\newcommand{\Ga}{\Gamma}

\newcommand{\ga}{\gamma}

\begin{document}

\title{\bf Toeplitz determinants from compatibility conditions}
\author{Estelle L. Basor \thanks{The first author was supported in part  
by National Science Foundation grant DMS-0200167 and also in part by  
the EPSRC for a Visiting Fellowship.}\\
          Department of Mathematics\\
          California Polytechnic University\\
          San Luis Obispo, CA 53407, USA\\
          ebasor@calpoly.edu
          \and
          Yang Chen\\
         Department of Mathematics\\
          Imperial College\\
          180 Queen's Gate\\
          London SW7 2BZ, UK\\
          ychen@ic.ac.uk}
\date{}

\maketitle

\begin{abstract}
In this paper we show, how a straightforward and natural
application of a pair of fundamental
identities valid for polynomials orthogonal over the
unit circle, can be used to
calculate the determinant of the finite Toeplitz matrix,
$$
\Delta_n=\det(w_{j-k})_{j,k=0}^{n-1}:=
\det\left(\int_{|z|=1}\frac{w(z)}{z^{j-k+1}}\frac{dz}{2\pi i}\right)
_{j,k=0}^{n-1},
$$
with the Fisher-Hartwig symbol,
$$
w(z)=C(1-z)^{\al+i\bt}(1-1/z)^{\al-i\bt},\quad |z|=1,\;\;\al>-1/2,
\;\bt\in{\mathbb R}\;.
$$
Here $C$ is the normalisation constant
chosen so that $w_0=\frac{1}{2\pi}.$ We use the same approach to
compute a difference equation for expressions related to the
determinants of the symbol 
$$w(z) = {\rm e}^{t(z+1/z)},$$ a symbol important
in the study of random permutations. Finally, we study the analogous  
equations for the symbol 
$$w(z) =  
{\rm e}^{tz}\prod_{\al=1}^{M}\left(\frac{z-a_{\al}}{z}\right)^{g_{\al}}.$$

\end{abstract}

\vfill\eject
\vskip .4cm
\setcounter{equation}{0}
\section{Introduction}
The large $n$ behavior of the determinant of $n\times n$ Toeplitz
matrices have seen many diverse physical applications,
from the Ising model \cite{MTW}, Random Matrix Theory \cite{Me},
String Theory \cite{GW,SP} and Combinatorics \cite{Ge,J, TW}. For an
arbitrary weight or symbol $w(z),$ the $n \times n$ Toeplitz matrix  is
defined by
\bq
T_{n}(w) = \left(\int_{|z|=1}\frac{w(z)}{z^{j-k+1}}\frac{dz}{2\pi
i}\right)
_{j,k=0}^{n-1}.
\eq
For smooth symbols $w$ the asymptotic behavior is a consequence of the
classis Szeg\"o Limit Theorem. For symbols that have jumps or other
kinds of singularities, the asymptotic behavior is harder to determine.
The earliest general conjecture for asymptotics of determinants for
such singular symbols dates back to the formulation of the the
Fisher-Hartwig conjecture. One of the key ingredients in the proof of
the conjecture was the exact computation of the determinants for
symbols of the form
\[
w_{\alpha, \beta}(z) = (1-z)^{\alpha + i\beta}(1 -1/z)^{\alpha
-i\beta}. \]
These symbols are called symbols with a pure Fisher-Hartwig
singularity. The exact computation was then combined with a
localization technique to find more general answers \cite{Basor}.
A good account of these results can be found in \cite{BS}. It is worth
noting that in the case of only jumps ($\alpha = 0$) the matrix reduces
to a Cauchy matrix and thus the determinants are straightforward to
calculate. The more general case was done by factoring the Toeplitz
matrix into a product of triangular and diagonal matrices \cite{Russian}.

Although, as mentioned above, asymptotic behavior in the case of smooth
symbols is given by the Strong Szeg\"o Limit Theorem, there are
applications where more exact information is desired. In the analysis
of problems involving random permutations certain difference equations
for quotients of the determinants also arise. These are related to the
symbol
$w(z) = e^{t(z + 1/z)}.$ For more results about the connection to
random permutations and random matrices see \cite{Ge, J, TW}.

The purpose of this paper is to re-derive these results and obtain new ones
using some very simple ideas from the theory of orthogonal polynomials. 
While the analysis of such singular symbols and smooth ones is generally very
different the two classes of examples share the property that the
derivative of their logarithms are rational functions. This is what
allows us to treat these classes as two examples of a general theory.

In the course of his investigation into the
inversion of finite Toeplitz matrices, Szeg\"o introduced orthogonal
polynomials supported on the unit circle and the corresponding
Szeg\"o kernel. To fix notations for the rest of the paper,
let
\bq
\phi_n(z)= k_nz^n+l_nz^{n-1}+...+\phi_n(0),\;\; k_n>0\;\;n=1,2,...
\eq
satisfy the orthogonality condition
\bq
\int_{|z|=1}\phi_{m}(z)\ov{\phi_{n}(z)}w(z)\frac{dz}{iz}=\de_{m,n},\;\;
m,n=0,1,2,...
\eq
Note that with the normalisation on the weight, $k_0=\phi_0(0)=1.$
  From the orthogonality condition, we find, \cite{Sze}
\bq
\mbox{}
k_nz\phi_n(z)&=&k_{n+1}\phi_{n+1}(z)-\phi_{n+1}(0)\phi_{n+1}^*(z)\\
\mbox{} k_n\phi_{n+1}(z)&=&k_{n+1}z\phi_n(z)+\phi_{n+1}(0)\phi_n^*(z),
\eq
and
\bq
k_n\phi_n(0)\phi_{n+1}(z)+k_{n-1}\phi_{n+1}(0)\phi_{n-1}(z)
=\left(k_n\phi_{n+1}(0)+k_{n+1}\phi_n(0)z\right)\phi_{n}(z)
\eq
by eliminating $\phi_{n}^*(z)$ from (1.3) and (1.4). The $^*$ operation
is defined as follows: If
\bq
\pi(z):=\sum_{j=0}^{n}a_jz^j\nonumber
\eq
then
\bq
\pi^{*}(z):=\sum_{j=0}^{n}\ov{a_j}z^{n-j}.\nonumber
\eq
The coefficient of $z^{n-1},$ $l_n$  can be obtained from the first order
difference equation,
\bq
\frac{l_{n+1}}{k_{n+1}}=\frac{l_n}{k_n}+
\ov{\left(\frac{\phi_n(0)}{k_n}\right)}
\frac{\phi_{n+1}(0)}{k_{n+1}},
\eq
if we are able to find $k_n$ and $\phi_n(0)$ from the weight. A simple
calculation shows that the Toeplitz determinant is
\bq
\Delta_n[w]=\prod_{j=0}^{n-1}\frac{1}{2\pi k_j^2}.
\eq
There are, as one may expect, very few cases where explicit formulae
were found for the Toeplitz determinants. However, it is clear that in
order to find information about the Toeplitz determinants. a possible
approach is to find explicitly the $k_{j}s$.

Our first step in that attempt is the derivation of two fundamental
identities valid for all $z\in\cpo.$
We first introduce the functions
$A_n(z)$ and $B_n(z)$
given by
\bq
\mbox{} A_n(z)&=&n\frac{k_{n-1}}{k_n}+i\frac{k_{n-1}}{\phi_n(0)}z
\int_{|\xi|=1}\frac{v^{\pr}(z)-v^{\pr}(\xi)}{z-\xi}
\phi_n(\xi)\overline{\phi_{n}^{*}(\xi)}w(\xi)d\xi\nonumber\\
\mbox{} B_n(z)&=&\frac{k_n}{k_{n-1}}\frac{A_n(z)}{z}-\frac{n}{z}
-i\int_{|\xi|=1}\frac{v^{\pr}(z)-v^{\pr}(\xi)}{z-\xi}\phi_n(\xi)
{\overline{\phi_n(\xi)}}w(\xi)d\xi,\nonumber
\eq
where $w(z)=:{\rm exp}(-v(z)).$
Then we claim
\bq
\quad\quad\quad\quad\quad B_n(z)+B_{n-1}(z)
=\frac{k_{n-1}}{k_{n-2}}\frac{A_{n-1}(z)}{z}+
\frac{k_{n}}{k_{n-2}}\frac{\phi_{n-1}(0)}{\phi_n(0)}A_{n-1}(z)
-\frac{n-1}{z}-
v^{\prime}(z),\quad\quad(S_1)\nonumber
\eq
\bq
\mbox{} &&(B_{n+1}(z)-B_n(z))(zk_{n+1}\phi_n(0)+k_n\phi_{n+1}(0))\nonumber\\
\mbox{} &&\quad\quad =
k_n\phi_n(0)A_{n+1}(z)-\frac{k^2_{n-1}}{k_{n-2}}
\frac{\phi_{n-1}(0)\phi_{n+1}(0)}{\phi_n(0)}A_{n-1}(z)-k_{n+1}\phi_{n}(0).
\quad\;\;\;\;\;\;(S_2)\nonumber
\eq
The equations $(S_1)$ and $(S_2)$ are the circular analogues obtained
earlier \cite{Chen1,Chen2} in the case where polynomials are orthogonal
with respect to weights supported on the interval $[a,b].$ The functions
$A_n(z)$ and $B_n(z)$ appear in the differentiation formula
\bq
\left(\frac{d}{dz}+B_n(z)\right)\phi_n(z)=A_{n-1}(z)\phi_{n-1}(z)
\eq
\bq
\mbox{}\biggl(\frac{d}{dz}+B_{n-1}(z)-\frac{k_{n-1}}{k_{n-
2}}\frac{A_n(z)}{z}
&-&\frac{k_n}{k_{n-2}}\frac{\phi_{n-1}(0)}{\phi_n(0)}A_{n-1}(z)\biggr)
\phi_{n-1}(z)\nonumber\\
\mbox{} &=&\frac{k_n}{k_{n-2}}\frac{\phi_{n-1}(0)}{\phi_{n}(0)}
\frac{A_{n-1}(z)}{z}\phi_n(z),
\eq
which can be thought of as generalised ``creation'' and ``annihilation''
operators. Equations (1.9) and (1.10) can be derived \cite{IW}
in a manner completely analogous that of the interval case.

The equation $(S_1)$ is simply found by the orthogonality condition and
the circular analogue of the Christoffel-Darboux formula.
Let
$$
\Psi_n(z):=\left(\matrix{\phi_{n}(z)\cr
                              \phi_{n-1}(z)}\right).
$$
Regarding $(S_2),$ this is found by re-writing (1.9) and (1.10), and the
recurrence relations as
\bq
\mbox{} \Psi_n^{\pr}(z)&=&M_n(z)\Psi_n(z)\nonumber\\
\mbox{} \Psi_{n+1}(z)&=&U_{n}(z)\Psi_n(z),\nonumber
\eq
and demanding that above over-determined systems are compatible;
entirely
analogues to what was done in \cite{Chen2}. It is now clear that
if $w^{\pr}/w$ is a rational function then, $(S_1)$ and $(S_2)$ will
supply
the basic equations for the determination of $k_n$ and $\phi_n(0).$
Furthermore, by eliminating $\phi_{n-1}(z)$ from (1.9) and (1.10), we
arrive at a second order differential equation satisfy by $\phi_n(z):$
\bq
Y^{\pr\pr}(z,n)+P(z,n)Y^{\pr}(z,n)+Q(z,n)Y(z,n)=0,
\eq
where,
\bq
\mbox{} P(z,n)&=&-\frac{n-1}{z}-v^{\pr}(z)-\frac{A_n^{\pr}(z)}{A_n(z)}\\
\mbox{} Q(z,n)&=&B_n^{\pr}(z)-B_n(z)\frac{A_n^{\pr}(z)}{A_n(z)}+
B_n(z)B_{n-1}(z)-\frac{k_{n-1}}{k_{n-2}}\frac{A_{n-1}(z)}{z}B_n(z)\nonumber\\
\mbox{}&-&\frac{k_n}{k_{n-2}}
\frac{\phi_{n-1}(0)}{\phi_n(0)}A_{n-1}(z)B_n(z)
+\frac{k_{n-1}}{k_{n-2}}
\frac{\phi_{n-1}(0)}{\phi_n(0)}\frac{A_{n-1}(z)A_n(z)}{z}.
\eq
We now give one more version of $(S_{1})$ and $(S_{2})$ which will
also prove to be useful in what follows. It turns out to be profitable to 
use the alternative parametrisation,
\bq
\mbox{} r_n&=&\frac{\phi_n(0)}{k_n}\nonumber\\
\mbox{} m_n&=&\frac{k_n}{k_{n+1}},\nonumber\\
\mbox{} s_n&=&\frac{r_{n+1}}{r_n},\nonumber
\eq
and the equations $(S_1)$ and $(S_2)$ become,
$$
B_{n+1}(z)+B_n(z)=\frac{1}{m_{n-1}}\frac{A_n(z)}{z}+\frac{A_n(z)}{m_{n-1}s_n}
-\frac{n}{z}-v^{\pr}(z).\eqno(T_1)
$$
$$
(B_{n+1}(z)-B_n(z))(z+s_n)=m_nA_{n+1}(z)-
\frac{s_n}{s_{n-1}}\frac{m_{n-1}^2}{m_{n-2}}A_{n-1}(z)-1.
\eqno(T_2)
$$

%

%


\section{The Pure Fisher-Hartwig Symbol}
\setcounter{equation}{0}
Let us return to the pure Fisher-Hartwig symbol,
\[
w_{\alpha, \beta}(z) = (1-z)^{\alpha + i\beta}(1- 1/z)^{\alpha -  
i\beta}\]
where the arguments are chosen so that the first factor is
analytic in the interior of the disc and one at $z = 0$ and that the
second factor is analytic outside the circle and one at infinity. The
Toeplitz determinants for this symbol have been previously computed
exactly. This was done using a hypergeometric approach and a
factorization of the corresponding
Toeplitz operator discovered by Roch. See Theorem 6.20 in \cite{BS} for  
the factorization. An alternative direct approach was
also done by Widom many years ago, but only recently in printed form \cite{BW}. 
We now give a difference equation derivation.

For our purposes we normalize $w$ so that $w_{0} = 1/2\pi.$
Note also that $v'$ is rational and that
\bq
\mbox{} v^{\prime}(z)&=&-\frac{2\al}{z-1}+\frac{\al-i\bt}{z}\nonumber\\
\mbox{} \frac{v^{\pr}(z)-v^{\pr}(\xi)}{z-\xi}&=&\frac{2\al
}{(z-1)(\xi-1)}
+\frac{i\bt-\al}{z\xi},\nonumber
\eq
and thus
\bq
\mbox{}
A_n(z)&=&n\frac{k_{n-1}}{k_n}+
\frac{k_{n-1}}{\phi_n(0)}\left(2\al J_n-(\al+i\bt)I_n\right)
+\frac{2\al\frac{k_{n-1}}{\phi_n(0)}(J_n-I_n)}{z-1}\\
\mbox{} B_{n}(z)&=&\frac{i\bt-\al-n}{z}+\frac{2\al(1-L_n)}{z-1}
+\frac{k_{n}}{k_{n-1}}\frac{A_n(z)}{z},
\eq
where
\bq
\mbox{} I_n&:=&\int_{|\xi|=1}\phi_n(\xi)\overline{\phi_n^*(\xi)}
\frac{w(\xi)}{i\xi}
d\xi\nonumber\\
\mbox{}
J_n&:=&\int_{|\xi|=1}\frac{\phi_n(\xi)\overline{\phi_n^{*}(\xi)}}{1-\xi}
\frac{w(\xi)}{i\xi}d\xi\nonumber\\
\mbox{}
L_n&:=&\int_{|\xi|=1}\frac{\phi_n(\xi)\overline{\phi_n(\xi)}}{1-\xi}
\frac{w(\xi)}{i\xi}d\xi.\nonumber
\eq
Now we have 5 unknowns; $I_n,$ $J_n,$ $L_n,$ $k_n$ and $\phi_n(0).$
Notice that $A_n(z)$ is of the form $c_1+c_2/(z-1)$ and
$B_n(z)$ is of the form $d_1/z+d_2/(z-1).$ We now show how five quite
simple steps yield the difference equation for the $k_{n}s.$

\vskip .2cm
\noindent
{\bf Step 1. Take the limit $z\to\infty$ in $(T_1).$}

Since $B_{n}(z)$ has no constant term the constant term of $A_{n}$ must
be zero. Thus we immediately know that
\[A_{n}(z) = \frac{2\alpha m_{n-1}(J_{n}-I_{n})}{r_{n}(z-1)} \]

\vskip .2cm
\noindent
{\bf Step 2. Compare the coefficients of $1/z$ in $(T_1).$}

A simple computation yields
\[ i\beta - n -\alpha = \frac{2 \alpha(J_{n}-I_{n})}{r_{n}} \]
and this implies
\[A_{n}(z) = \frac{m_{n-1}(-\alpha +i\beta -n)}{(z-1)} .\]
This of course also gives information about $B_{n}$ and at this  point
we can conclude that
\[ B_{n}(z) = \frac{2\alpha(1-L_{n}) -\alpha +i\beta -n}{z-1}.
\]
\vskip .2cm
\noindent
{\bf Step 3. Compare the coefficients of $1/(z-1)$ in $(T_1).$}
 
From this we have
\bq
4\al- 2\al (L_{n+1}+ L_{n})  -\alpha +i\beta -n -1=
\frac{r_{n}}{r_{n+1}}(-\alpha +i\beta -n)+2\al.
\eq
We will return to this equation in a moment.
\vskip .2cm
\noindent
{\bf Step 4. Take the limit $z\to\infty$ in $(T_2).$}

This limit shows
\[-1 -2\alpha(L_{n+1}-L_{n})= -1\]
or that $L_{n+1} = L_{n}.$ A  direct computation shows that $L_{0} =
\frac{\alpha +i\beta}{2\alpha}.$
Thus $B_{n}(z) = \frac{-n}{z-1}$ and returning to  step three we see
that we have
an equation for $r_{n},$
namely,
\[-\alpha -i\beta  -n -1= \frac{r_{n}}{r_{n+1}}( -\alpha + i\beta  -n
).\]
\vskip .2cm
\noindent
{\bf Step 5. Compare the coefficients of $1/(z-1)$ in $(T_2).$}

We have that
$s_{n} = \frac{r_{n+1}}{r_{n}}$ so this residue produces the equation
\[ 1+\frac{r_{n+1}}{r_{n}} = m_{n}^{2}(\alpha -i\beta + n +1) -
\frac{r_{n+1}}{r_{n}}\frac{r_{n-1}}{r_{n}}m_{n-1}^{2}(\alpha +i\beta +n
-1) \]
and since we have an expression for the quotients of the $r_{n}s$ this
becomes
\bq
\mbox{}
2\al+2n+1&=&m_n^2(\al+i\bt+n+1)(\al-i\bt+n+1)-m_{n-1}^2(\al+i\bt+n)
(\al-i\bt+n),\nonumber\\
\eq
We easily verify that,
\bq
m_n^2=\frac{(n+1)(2\al+n+1)}{(\al+i\bt+n+1)(\al-i\bt+n+1)},
\eq
solves (2.9) with the initial condistions $k_{-1}=0$ and $k_0=1.$
 From this it follows
\bq
\mbox{} k_n^2&=&\frac{\Ga(2\al+1)}{\Ga(\al-i\bt+1)\Ga(\al+i\bt+1)}
\frac{\Ga(\al+i\bt+n+1)}{\Ga(n+1)}\frac{\Ga(\al-
i\bt+n+1)}{\Ga(n+2\al+1)}\\
\mbox{} |\phi_n(0)|^2&=&\frac{\Ga(2\al+1)}{\Ga(n+1)\Ga(n+2\al+1)}
\frac{|\Ga(\al+i\bt+n)|^2}{|\Ga(\al+i\bt)|^2}\\
\mbox{} l_n&=&\frac{(\al-i\bt)n}{n+\al+i\bt}k_n.
\eq
Incidentally $I_n$ and $J_n$ can also be easily determined. In Step 1,
the vanishing of the constant term in $A_n(z)$ implies
\bq
n+\frac{2\al J_n-(\al+i\bt)I_n}{r_n}=0,\nonumber
\eq
and when combined with the first equation in Step 2 gives 
\bq
\mbox{} I_n&=&r_n\nonumber\\
\mbox{} J_n&=&\frac{(\al+i\bt-n)r_n}{2\al}.\nonumber
\eq

\setcounter{equation}{0}
{\section {Toeplitz determinant and discriminant.}}

The computation of the Toeplitz determinant $\Delta_n$, is now immediate
\bq
\mbox{} \Delta_n&=&\prod_{j=0}^{n-1}\frac{1}{2\pi k_j^2}\nonumber\\
\mbox{} &=&C^n\frac{G(n+1)G(n+2\al+1)G(\al+i\bt+1)G(\al-i\bt+1)}
{G(2\al+1)G(n+\al+i\bt+1)G(n+\al-i\bt+1)},
\eq
where $G(z)$ is Barnes $G-$function.
Let
\bq
\pi_n(z)=\ga\prod_{j=1}^{n}(z-z_j(n)),\nonumber
\eq
where $z_j(n)$ are the $n$ (simple) zeros, then the discriminant is
\bq
D[\pi_n]=\ga^{2n-2}\prod_{1\leq j<k\leq n}(z_j(n)-z_k(n))^2.
\eq
Indeed (see (5.15) of \cite{IW})
\bq
\mbox{}
D[\phi_n]&=&(-1)^{n(n-1)/2}\frac{(\phi_n(0))^{n-1}}{k_nk_{n-1}^n}
\;\prod_{j=1}^{n-1}k_j^2\;\;\prod_{l=1}^{n}A_n(z_l(n))\nonumber\\
\mbox{} &=&(-1)^{n(n-1)/2}
\left(\frac{\phi_n(0)}{k_n}\right)^{n-1}\frac{1}{k_n^2}
\left(\frac{k_n}{k_{n-1}}\right)^n\;\;\prod_{j=1}^{n-1}k_j^2
\;\;\prod_{l=1}^{n}A_n(z_l(n))\nonumber\\
\mbox{} &=&(-1)^{n(n-1)/2}\frac{r_n^{n-1}}{m_{n-1}^nk_n^2}
\;\;\prod_{j=1}^{n-1}k_j^2\;\;\prod_{l=1}^{n}A_n(z_l(n)).
\eq
The only missing information required for the computation of $D[\phi_n]$ is
$$
\prod_{l=1}^{n}A_n(z_l).
$$
The derivative of $\phi_n(z)$ is (see (1.10)),
\bq
\phi_n^{\pr}(z)=A_n(z)\phi_{n-1}(z)-B_n(z)\phi_n(z),
\eq
where in our example,
\bq
\mbox{} A_n(z)&=&\frac{a_n}{z-1},\;\;a_n=m_{n-1}(i\bt-\al-n)\nonumber\\
\mbox{} B_n(z)&=&\frac{b_n}{z-1},\;\;b_n=-n.\nonumber
\eq
Now
\bq
\prod_{l=1}^{n}A_n(z_l)=a_n^n\prod_{l=1}^{n}\frac{1}{z_l-1}
=a_n^n(-1)^n\frac{k_n}{\phi_n(1)},\nonumber
\eq
from the fact that $\phi_n(z)=k_n\prod_{l=1}^{n}(z-z_l).$ To determine
$\phi_n(1),$ evaluate both sides of (3.4) at $z=1;$ keeping in mind that
the l.h.s. is regular at every finite $z.$ Therefore the residue of the 
r.h.s. at $z=1$ must vanish:
\bq
\mbox{} \phi_n(1)&=&\frac{a_n}{b_n}\phi_{n-1}(1)\nonumber\\
\mbox{} &=&\prod_{j=1}^{n}\frac{a_j}{b_j},\nonumber\\
\mbox{}
\prod_{l=1}^{n}A_n(z_l)&=&(-1)^na_n^n\frac{k_n}{\phi_n(1)},\nonumber\\
\mbox{} D[\phi_n]&=&(-1)^{n(n+1)/2}(i\bt-\al-n)^n
\frac{r_n^{n-1}}{k_n}\left(\prod_{j=1}^{n-1}k_j^2\right)
\left(\prod_{j=1}^{n}\frac{b_j}{a_j}\right)\nonumber\\
\mbox{} &=&(-1)^{n(n+3)/2}\left(\frac{n+\al-i\bt}{2\pi}\right)^n
\frac{r_n^{n-1}}{k_n}\;\Delta_n^{-1}
\left(\prod_{j=1}^{n}\frac{b_j}{a_j}\right)
\eq
Now
\bq
\mbox{} \prod_{j=1}^{n}\frac{b_j}{a_j}&=&
{\sqrt {\prod_{j=1}^{n}\frac{\al+i\bt+j}{\al-i\bt+j}\;\frac{j}{j+2\al}}}
\nonumber\\
\mbox{} &=&{\sqrt {\frac{\Ga(n+1+\al+i\bt)}{\Ga(n+1+\al-i\bt)}
\frac{\Ga(\al-i\bt+1)}{\Ga(\al+i\bt+1)}\frac{\Ga(2\al+1)\Ga(n+1)}
{\Ga(n+2\al+1)}}}.
\eq
Therefore,
\bq
\Delta_nk_n|D[\phi_n]|=\left(\frac{|n+\al-i\bt|}{2\pi}\right)^n
|r_n|^{n-1}\;\Biggl|\prod_{j=1}^{n}\frac{b_j}{a_j}\Biggr|.
\eq
Note that according to general theory \cite{Sze},
$$\lim_{n\to\infty}k_n=\kappa>0,$$
and in particular $$\lim_{n\to\infty}\phi_n(0)=0.$$ Indeed this is the
case for the Fisher-Hartwig symbol.

   To determine
$|D[\phi_n]|,$ as $n\to\infty,$ we give here some preliminary results:
\bq
\mbox{}\lim_{n\to\infty}k_n&=&\frac{\sqrt
{\Ga(2\al+1)}}{|\Ga(\al+i\bt+1)|}
=:\kappa,
\nonumber\\
\mbox{} \phi_n(0)&\sim&\kappa\frac{\Ga(\al+i\bt+1)}{\Ga(\al-i\bt)}
n^{-1-2i\bt},\nonumber\\
\mbox{}
|r_n|^{n-1}&\sim&\left(\frac{|\al+i\bt|}{n}\right)^{n-1},\nonumber\\
\mbox{} \Biggl|\prod_{j=1}^{n}\frac{b_j}{a_j}\Biggr|
&\sim&\frac{{\sqrt {\Ga(2\al+1)}}}{n^{\al}},\nonumber\\
\mbox{} \Delta_n&\sim&
\frac{|G(\al+i\bt+1)|^2}{G(2\al+1)}C^{n}n^{\al^2+\bt^2}.
\nonumber
\eq
Therefore,
\bq
|D[\phi_n]|\sim\;\frac{|\Ga(\al+i\bt+1)|}{|\al+i\bt|}\;
\frac{G(2\al+1)}{|G(\al+i\bt+1)|^2}
\left(\frac{\Ga(2\al+1)}{C|\al+i\bt||\Ga(\al+i\bt)|^2}\right)^n\;n^{1-
\al-\al^2-\bt^2}.
\eq

\setcounter{equation}{0}
{\section {Differential Equation.}}
  From the general theory $\phi_n(z)$ satisfies
\bq
Y^{\pr\pr}(z)+P(z,n)Y^{\pr}(z)+Q(z,n)Y(z)=0,\nonumber
\eq
where
\bq
\mbox{} P(z,n)&=&\frac{1-n-\al+i\bt}{z}+\frac{2\al+1}{z-1},\nonumber\\
\mbox{} Q(z,n)&=&-\frac{n(\al+i\bt+1)}{z(z-1)}.
\mbox{}
\eq
The general solution of the differential equation is
\bq
\mbox{} &A&\;_2F_1(-n,\al+i\bt+1;1-n-\al+i\bt;z)\nonumber\\
\mbox{} +&B&\;z^{n+\al-i\bt}\;_2F_1(n+2\al+1,\al-i\bt;n+\al-i\bt+1;z),
\eq
and since we know that $\phi_{n}$ is a polynomial we have that
\[
\phi_{n}(z) = A\,_{2}F_1(-n,\al+i\bt+1;1-n-\al+i\bt;z)
\]
where 
\[
A  = \sqrt{\frac{\Gamma(2\al+1)\Gamma(\al+i\beta +n+1)\Gamma(\al-i\beta +n+1)}
{\Gamma(\al-i\beta +a)\Gamma(\al+i\beta +1)\Gamma(n+1)\Gamma(n+2\al+1)}}
\frac{\Gamma(n+\al-i\beta)\Gamma(\al+i\beta+1)}{\Gamma(n+\al+i\beta+1)\Gamma(\al-i\beta)}. 
\]
This value comes from computing the coefficient for the $nth$ term of the Hypergeometric function and also using the value of  $k_{n}$  given in (2.6). 
Indeed this particular Hypergeometric function was found by Askey to be
orthogonal with respect to the pure Fisher-Hartwig symbol in a
commentary on Szeg\"o's collected papers \cite{Askey}.  

\setcounter{equation}{0}
{\section {The weight $C{\rm exp}[t(z+1/z)/2]$.}}

For this example
\bq
&&v^{\pr}(z)=-\frac{t}{2}+\frac{t}{2z^2}\nonumber\\
&&\frac{v^{\pr}(z)-v^{\pr}(\xi)}{z-\xi}=-\frac{1}{z}\frac{t}{2\xi^2}
-\frac{1}{z^2}\frac{t}{2\xi}.\nonumber
\eq
We find,
\bq
A_n(z)
&=&(n+b_n)m_{n-1}+\frac{m_{n-1}a_n}{z}.
\nonumber\\
B_n(z)
&=&\frac{L_n+b_n}{z}+\frac{a_n-t/2}{z^2}\nonumber\\
a_n&=&\frac{t}{2r_n}
\int_{|\xi|=1}\frac{\phi_n(\xi)\ov{\phi_n^*(\xi)}}{i\xi}w(\xi)
d\xi\nonumber\\
b_n&=&\frac{t}{2r_n}\int_{|\xi|=1}\frac{\phi_n(\xi)\ov{\phi_n^*(\xi)}}
{i\xi^2}w(\xi)d\xi\nonumber\\
L_n&=&-
\frac{t}{2}\int_{|\xi|=1}\frac{\phi_n(\xi)\ov{\phi_n(\xi)}}{i\xi^2}
w(\xi)d\xi\nonumber.
\eq
There are 5 unknowns: $r_n$, $m_n,$ $a_n,$ $b_n,$ and $L_n.$ Again we
use the same basic steps as was done in the previous example.
\vskip .2cm
\noindent
{\bf Step 1. Take the limit as $z\to\infty$ in $(T_1).$}
\bq
b_n=-n-\frac{t}{2}s_n.
\eq
{\bf Step 2. Compare the coefficients of $1/z$ in $(T_1)$.}
\bq
L_{n+1}+L_n-\frac{t}{2}\left(s_{n+1}+\frac{2a_n/t}{s_n}\right)=n+1.
\eq
{\bf Step 3. Compare the coefficients of $1/z^2$ in $(T_1)$.}
\bq
a_{n+1}+a_n-t&=&a_n-t/2\nonumber\\
a_n&=&\frac{t}{2}.
\eq
{\bf Step 4. Take the limit as $z\to\infty$ in $(T_2)$.}
This limit yields after simplifying
\bq
L_{n+1}-\frac{t}{2}s_{n+1}(1-m_n^2)-\left(L_n-\frac{t}{2}s_n(1-m_{n-
1}^2)
\right)=0\nonumber
\eq
which implies,
\bq
L_n=\al+\frac{t}{2}s_n(1-m^2_{n-1}),\nonumber
\eq
where $\al$ is an ``integration'' constant. To determine $\al$ put $n=0$
and note that $m_{-1}=0$ and $s_0=r_1/r_0=r_1=-I_1(t)/I_0(t),$ where
$I_j(t)$
is the $I-$Bessel function of order $j.$ So,
\bq
\al=L_0-\frac{t}{2}r_1.\nonumber
\eq
But
\bq
L_0&=&-\frac{tC}{2}\int_{|z|=1}\frac{\exp(t(z+1/z)/
2)}{iz^2}dz\nonumber\\
&=&Ct\int_{0}^{\pi}{\rm e}^{-t\cos\psi}\cos\psi d\psi\nonumber\\
&=&-\pi CtI_1(t)=-\frac{\pi tI_1(t)}{2\pi I_0(t)}=-\frac{t}{2}
\frac{I_1(t)}{I_0(t)}.\nonumber
\eq
So $\al=0,$ and
\bq
L_n=\frac{t}{2}s_n(1-m^2_{n-1}).
\eq
{\bf Step 5. Compare the coefficients of $1/z$ in $(T_2)$.}
Again after simplifying
\bq
L_{n+1}-\frac{t}{2}\left(s_{n+1}+\frac{m^2_n}{s_n}\right)
-\left(L_{n}-\frac{t}{2}\left(s_n+\frac{m^2_{n-1}}{s_{n-
1}}\right)\right)=1,
\nonumber
\eq
which implies,
\bq
L_n-\frac{t}{2}\left(s_n+\frac{m^2_{n-1}}{s_{n-1}}\right)=\bt+n.
\nonumber
\eq
It turns out that the ``integration'' constant $\bt$ is also 0.
Therefore,
\bq
L_n=n+\frac{t}{2}\left(s_n+\frac{m^2_{n-1}}{s_{n-1}}\right).
\eq
{\bf Step 6. Compare coefficients of $1/z^2$ in $(T_2)$}
This gives $0=0.$

Eliminating $L_n$ from (5.4) and (5.5) gives,
\bq
-\frac{2n}{t}=m^2_{n-1}\left(s_n+\frac{1}{s_{n-1}}\right).
\eq
Now substitute (5.5), (5.4) into (5.2) and note that $a_n=t/2,$ we get
after some simplification,
\bq
s_n^2=\frac{1-m^2_n}{1-m^2_{n-1}}.
\eq
Now since $s_n=r_{n+1}/r_n,$ we see that
\bq
\frac{r_{n+1}^2}{1-m_n^2}=\frac{r_n^2}{1-m_{n-1}^2}=\gamma={\rm constant.}
\nonumber
\eq
which implies,
\bq
m^2_{n-1}=1-r^2_n,
\eq
where an easy computation shows that $\gamma=1.$
This last equation hold in general \cite{Ger}.

Using this on (5.6) gives the discrete Painleve II \cite{Hi,SP}
\bq
r_{n+1}+r_{n-1}=-\frac{2n}{t}\;\frac{r_n}{1-r_n^2}.
\eq
\setcounter{equation}{0}
{\section {Another non-linear difference equation for $r_n$.}}

In \cite{ITW} the authors considered a symbol of the form
\bq
w(z)={\rm e}^{tz}\prod_{\al=1}^{M}\left(\frac{z-z_{\al}}{z}
\right)^{g_{\al}},\quad \sum_{\al=1}^{M}g_\al=g>0,\;\;\;-1<z_{\al}<0,
\;\;\;\sum_{\al=1}^{M}g_\al z_\al=-1.
\eq
We have,
\bq
v(z)&=&-tz-\sum_{\al}g_\al\left(\ln(z-z_\al)-\ln z\right)\nonumber\\
v^{\pr}(z)&=&-t-\sum_{\al}\frac{g_\al}{z-z_\al}+\frac{g}{z}\nonumber\\
\frac{v^{\pr}(z)-v^{\pr}(\xi)}{z-\xi}&=&-\frac{g}{z\xi}+\sum_{\al}
\frac{g_\al}{z-z_\al}\:\frac{1}{\xi-z_\al}.\nonumber
\eq
\bq
A_n(z)&=&(n+a_{n})m_{n-1}+m_{n-1}\sum_\al\frac{b_n(\al)}{z-
z_\al}\nonumber\\
a_n&=&\frac{1}{r_n}\int_{|\xi|=1}
\left(\frac{g}{\xi}-\sum_\al\frac{g_\al}{\xi-z_\al}\right)\phi_n(\xi)
\ov{\phi^*_n(\xi)}\frac{w(\xi)}{i}d\xi\nonumber\\
b_n(\al)&=&-\frac{g_\al z_\al}{r_n}\int_{|\xi|=1}
\frac{\phi_n(\xi)\ov{\phi^*_n(\xi)}}{\xi-z_\al}\frac{w(\xi)}{i}d\xi.
\nonumber\\
B_{n}(z)&=&\frac{A_n(z)}{m_{n-1}z}-\frac{n+g}{z}+\sum_{\al}
\frac{g_\al+c_n(\al)}{z-z_\al}\nonumber\\
c_n(\al)&=&g_\al z_\al\int_{|\xi|=1}\frac{\phi_n(\xi)\ov{\phi_n(\xi)}}
{\xi-z_\al}\frac{w(\xi)}{i\xi}d\xi.
\nonumber\eq
{\bf Step 1.} $z\to\infty$ in $(T_1):$
\bq
a_n+n=-\frac{r_{n+1}}{r_n}t=-s_nt.
\eq
{\bf Step 2.} Residues at $z=0$ of $(T_1):$
\bq
a_n=g+\sum_{\al}\frac{b_n(\al)}{z_\al}.
\eq
{\bf Step 3.} Residues at $z=z_\al$ of $(T_1):$
\bq
g_\al+c_{n+1}(\al)+c_n(\al)+\frac{b_{n+1}(\al)}{z_\al}=\frac{b_n(\al)}{s
_n}.
\eq
{\bf Step 4.} $z\to\infty$ in $(T_2):$
\bq
1+a_{n+1}-a_n+\sum_{\al}\left(c_{n+1}(\al)-c_n(\al)\right)
=m_n^2(n+1+a_{n+1})-\frac{s_n}{s_{n-1}}m^2_{n-1}(n-1+a_{n-1})\nonumber
\eq
which is simplifies to:
\bq
tr^2_ns_n-tr^2_{n+1}+\sum_{\al}\left(c_{n+1}(\al)-
c_n(\al)\right)=0\nonumber
\eq
using (6.2) and $m^2_{n-1}=1-r_n^2.$ ``Integrating'' the above gives:
\bq
tr_{n}r_{n+1}-\sum_{\al}c_n(\al)=\la(\al),
\eq
where $\la(\al)$ is the $n-$ independent integration constant and  
depends on $t$ and the set of $z_{\al}'s$.

\noindent
{\bf Step 5.} Residues at $z=0$ of $(T_2):$
\bq
a_{n+1}-a_{n}=\sum_{\al}\frac{b_{n+1}(\al)-b_{n}(\al)}{z_\al},\nonumber
\eq
and upon ``integration'' gives,
\bq
a_n-\sum_{\al}\frac{b_n(\al)}{z_\al}=d(\al)=g,
\eq
the same as (6.3).

\noindent
{\bf Step 6.} Residues at $z=z_\al$ of $(T_2):$
\bq
(z_\al+s_n)\left(c_{n+1}(\al)-c_{n}(\al)+\frac{b_{n+1}(\al)-
b_n(\al)}{z_\al}\right)
=m^2_nb_{n+1}(\al)-\frac{s_n}{s_{n-1}}m^2_{n-1}b_{n-1}(\al).
\eq

Let's see what all this means when there is just one $z_{\al}$.

If we use Step 3, along with (6.5) and (6.6) we have that

\[
g  + tr_{n+2}r_{n+1} + \lambda + tr_{n}r_{n+1} +\lambda -g - s_{n+1}t  
-n-1 =  \frac{z_{\al}}{s_{n}}(-g
-s_{n}t -n),
\]
or, after simplifying,
\[
t r_{n+1}r_{n} + + tr_{n}r_{n+1} +2\lambda  - s_{n+1}t -n-1 =  1/s_{n}  
-tz_{\al}-nz_{\al}/s_{n},
\]
A little more algebra yields
\[
r_{n+2} + r_{n} = -
\frac{(nr_{n} -t)z_{\al} -2\lambda + (n+1)}{t(1- r_{n+1}^{2})}
\]
an equation very much like the one found for the previous example.

\end{document}